\documentclass[final,nomarks]{dmtcs-episciences}

\usepackage[utf8]{inputenc}
\usepackage{subfigure}

\usepackage{amsmath,amsthm}
\usepackage{amssymb}
\usepackage{color}
\usepackage{xspace}
\def\red{\textcolor{red} }

\def\c{\ensuremath{\mathcal C}\xspace}

\def\i{\ensuremath{\mathcal I}\xspace}
\def\eps{\epsilon}
\def\k{blue\xspace}

\def\lr{LRMax\xspace}

\def\gl{ground level\xspace}

\def\avs{$\{4321,\,3241\}$-avoiders\xspace}
\def\gf{generating function\xspace}

\def\mbf#1{\mathchoice{\hbox{\boldmath $\displaystyle #1$}}
        {\hbox{\boldmath $\textstyle #1$}}
        {\hbox{\boldmath $\scriptstyle #1$}}
        {\hbox{\boldmath $\scriptscriptstyle #1$}}} 

\newcommand{\cc}[2]{C_{#1}^{(#2)}}

\catcode`\@=11

\thicklines
\newskip\Einheit \Einheit=.6cm
\newcount\xcoord \newcount\ycoord
\newdimen\xdim \newdimen\ydim \newdimen\PfadD@cke \newdimen\Pfadd@cke
\def\PfadDicke#1{\PfadD@cke#1 \divide\PfadD@cke by2 
\Pfadd@cke\PfadD@cke \multiply\PfadD@cke by2}
\long\def\LOOP#1\REPEAT{\def\BODY{#1}\ITERATE}
\def\ITERATE{\BODY \let\next\ITERATE \else\let\next\relax\fi \next}
\let\REPEAT=\fi
\def\Punkt{\hbox{\raise-2pt\hbox to0pt{\hss\scriptsize$\bullet$\hss}}}

\def\DuennPunkt(#1,#2){\unskip
  \raise#2 \Einheit\hbox to0pt{\hskip#1 \Einheit
          \raise-1.5pt\hbox to0pt{\hss\tiny$\bullet$\hss}\hss}}
		  
\def\NormalPunkt(#1,#2){\unskip
  \raise#2 \Einheit\hbox to0pt{\hskip#1 \Einheit
          \raise-3pt\hbox to0pt{\hss\large$\bullet$\hss}\hss}}
\def\DickPunkt(#1,#2){\unskip
  \raise#2 \Einheit\hbox to0pt{\hskip#1 \Einheit
          \raise-4pt\hbox to0pt{\hss\Large$\bullet$\hss}\hss}}
\def\Kreis(#1,#2){\unskip
  \raise#2 \Einheit\hbox to0pt{\hskip#1 \Einheit
          \raise-4pt\hbox to0pt{\hss\Large$\circ$\hss}\hss}}
\def\Diagonale(#1,#2)#3{\unskip\leavevmode
  \xcoord#1\relax \ycoord#2\relax
      \raise\ycoord \Einheit\hbox to0pt{\hskip\xcoord \Einheit
         \unitlength\Einheit
         \line(1,1){#3}\hss}}
\def\AntiDiagonale(#1,#2)#3{\unskip\leavevmode
  \xcoord#1\relax \ycoord#2\relax \advance\xcoord by -0.05\relax
      \raise\ycoord \Einheit\hbox to0pt{\hskip\xcoord \Einheit
         \unitlength\Einheit
         \line(1,-1){#3}\hss}}
\def\Pfad(#1,#2),#3\endPfad{\unskip\leavevmode
  \xcoord#1 \ycoord#2 \thicklines\ZeichnePfad#3\endPfad\thinlines}
\def\ZeichnePfad#1{\ifx#1\endPfad\let\next\relax
  \else\let\next\ZeichnePfad
    \ifnum#1=1
      \raise\ycoord \Einheit\hbox to0pt{\hskip\xcoord \Einheit
         \vrule height\Pfadd@cke width1 \Einheit depth\Pfadd@cke\hss}%
      \advance\xcoord by 1
     \else\ifnum#1=2
      \raise\ycoord \Einheit\hbox to0pt{\hskip\xcoord \Einheit
         \unitlength\Einheit
         \line(0,1){1}\hss}
      \advance\xcoord by 0
      \advance\ycoord by 1
 \else\ifnum#1=3
      \raise\ycoord \Einheit\hbox to0pt{\hskip\xcoord \Einheit
         \unitlength\Einheit
         \line(1,1){1}\hss}
      \advance\xcoord by 1
      \advance\ycoord by 1
    \else\ifnum#1=4
      \raise\ycoord \Einheit\hbox to0pt{\hskip\xcoord \Einheit
         \unitlength\Einheit
         \line(1,-1){1}\hss}
      \advance\xcoord by 1
      \advance\ycoord by -1
   \else\ifnum#1=5
      \raise\ycoord \Einheit\hbox to0pt{\hskip\xcoord \Einheit
         \unitlength\Einheit
         \line(2,1){2}\hss}
      \advance\xcoord by 2
      \advance\ycoord by 1
	  \else\ifnum#1=6
      \raise\ycoord \Einheit\hbox to0pt{\hskip\xcoord \Einheit
         \unitlength\Einheit
         \line(2,-1){2}\hss}
      \advance\xcoord by 2
      \advance\ycoord by -1
	  \else\ifnum#1=7
      \raise\ycoord \Einheit\hbox to0pt{\hskip\xcoord \Einheit
         \unitlength\Einheit
         \line(3,1){3}\hss}
      \advance\xcoord by 3
      \advance\ycoord by 1
	  \else\ifnum#1=8
      \raise\ycoord \Einheit\hbox to0pt{\hskip\xcoord \Einheit
         \unitlength\Einheit
         \line(3,-1){3}\hss}
      \advance\xcoord by 3
      \advance\ycoord by -1
    \fi\fi\fi\fi\fi\fi\fi\fi
  \fi\next}
\def\hSSchritt{\leavevmode\raise-.4pt\hbox 
to0pt{\hss.\hss}\hskip.2\Einheit
  \raise-.4pt\hbox to0pt{\hss.\hss}\hskip.2\Einheit
  \raise-.4pt\hbox to0pt{\hss.\hss}\hskip.2\Einheit
  \raise-.4pt\hbox to0pt{\hss.\hss}\hskip.2\Einheit
  \raise-.4pt\hbox to0pt{\hss.\hss}\hskip.2\Einheit}
\def\vSSchritt{\vbox{\baselineskip.2\Einheit\lineskiplimit0pt
\hbox{.}\hbox{.}\hbox{.}\hbox{.}\hbox{.}}}
\def\DSSchritt{\leavevmode\raise-.4pt\hbox to0pt{%
  \hbox to0pt{\hss.\hss}\hskip.2\Einheit
  \raise.2\Einheit\hbox to0pt{\hss.\hss}\hskip.2\Einheit
  \raise.4\Einheit\hbox to0pt{\hss.\hss}\hskip.2\Einheit
  \raise.6\Einheit\hbox to0pt{\hss.\hss}\hskip.2\Einheit
  \raise.8\Einheit\hbox to0pt{\hss.\hss}\hss}}
\def\dSSchritt{\leavevmode\raise-.4pt\hbox to0pt{%
  \hbox to0pt{\hss.\hss}\hskip.2\Einheit
  \raise-.2\Einheit\hbox to0pt{\hss.\hss}\hskip.2\Einheit
  \raise-.4\Einheit\hbox to0pt{\hss.\hss}\hskip.2\Einheit
  \raise-.6\Einheit\hbox to0pt{\hss.\hss}\hskip.2\Einheit
  \raise-.8\Einheit\hbox to0pt{\hss.\hss}\hss}}
\def\SPfad(#1,#2),#3\endSPfad{\unskip\leavevmode
  \xcoord#1 \ycoord#2 \ZeichneSPfad#3\endSPfad}
\def\ZeichneSPfad#1{\ifx#1\endSPfad\let\next\relax
  \else\let\next\ZeichneSPfad
    \ifnum#1=1
      \raise\ycoord \Einheit\hbox to0pt{\hskip\xcoord \Einheit
         \hSSchritt\hss}%
      \advance\xcoord by 1
    \else\ifnum#1=2
      \raise\ycoord \Einheit\hbox to0pt{\hskip\xcoord \Einheit
        \hbox{\hskip-2pt \vSSchritt}\hss}%
      \advance\ycoord by 1
    \else\ifnum#1=3
      \raise\ycoord \Einheit\hbox to0pt{\hskip\xcoord \Einheit
         \DSSchritt\hss}
      \advance\xcoord by 1
      \advance\ycoord by 1
    \else\ifnum#1=4
      \raise\ycoord \Einheit\hbox to0pt{\hskip\xcoord \Einheit
         \dSSchritt\hss}
      \advance\xcoord by 1
      \advance\ycoord by -1
    \fi\fi\fi\fi
  \fi\next}
\def\Koordinatenachsen(#1,#2){\unskip
 \hbox to0pt{\hskip-.5pt\vrule height#2 \Einheit width.5pt depth1 
\Einheit}%
 \hbox to0pt{\hskip-1 \Einheit \xcoord#1 \advance\xcoord by1
    \vrule height0.25pt width\xcoord \Einheit depth0.25pt\hss}}
\def\Koordinatenachsen(#1,#2)(#3,#4){\unskip
 \hbox to0pt{\hskip-.5pt \ycoord-#4 \advance\ycoord by1
    \vrule height#2 \Einheit width.5pt depth\ycoord \Einheit}%
 \hbox to0pt{\hskip-1 \Einheit \hskip#3\Einheit 
    \xcoord#1 \advance\xcoord by1 \advance\xcoord by-#3 
    \vrule height0.25pt width\xcoord \Einheit depth0.25pt\hss}}
\def\Gitter(#1,#2){\unskip \xcoord0 \ycoord0 \leavevmode
  \LOOP\ifnum\ycoord<#2
    \loop\ifnum\xcoord<#1
      \raise\ycoord \Einheit\hbox to0pt{\hskip\xcoord 
\Einheit\Punkt\hss}%
      \advance\xcoord by1
    \repeat
    \xcoord0
    \advance\ycoord by1
  \REPEAT}
\def\Gitter(#1,#2)(#3,#4){\unskip \xcoord#3 \ycoord#4 \leavevmode
  \LOOP\ifnum\ycoord<#2
    \loop\ifnum\xcoord<#1
      \raise\ycoord \Einheit\hbox to0pt{\hskip\xcoord 
\Einheit\Punkt\hss}%
      \advance\xcoord by1
    \repeat
    \xcoord#3
    \advance\ycoord by1
  \REPEAT}
\def\Label#1#2(#3,#4){\unskip \xdim#3 \Einheit \ydim#4 \Einheit
  \def\lo{\advance\xdim by-.5 \Einheit \advance\ydim by.5 \Einheit}%
  \def\llo{\advance\xdim by-.25cm \advance\ydim by.5 \Einheit}%
  \def\loo{\advance\xdim by-.5 \Einheit \advance\ydim by.25cm}%
  \def\o{\advance\ydim by.25cm}%
  \def\ro{\advance\xdim by.5 \Einheit \advance\ydim by.5 \Einheit}%
  \def\rro{\advance\xdim by.25cm \advance\ydim by.5 \Einheit}%
  \def\roo{\advance\xdim by.5 \Einheit \advance\ydim by.25cm}%
  \def\l{\advance\xdim by-.30cm}%
  \def\r{\advance\xdim by.30cm}%
  \def\lu{\advance\xdim by-.5 \Einheit \advance\ydim by-.6 \Einheit}%
  \def\llu{\advance\xdim by-.25cm \advance\ydim by-.6 \Einheit}%
  \def\luu{\advance\xdim by-.5 \Einheit \advance\ydim by-.30cm}%
  \def\u{\advance\ydim by-.30cm}%
  \def\ru{\advance\xdim by.5 \Einheit \advance\ydim by-.6 \Einheit}%
  \def\rru{\advance\xdim by.25cm \advance\ydim by-.6 \Einheit}%
  \def\ruu{\advance\xdim by.5 \Einheit \advance\ydim by-.30cm}%
  #1\raise\ydim\hbox to0pt{\hskip\xdim
     \vbox to0pt{\vss\hbox to0pt{\hss$#2$\hss}\vss}\hss}%
}
\catcode`\@=12

\newtheorem{theorem}{Theorem}

\newtheorem{lemma}[theorem]{Lemma}
\newtheorem{prop}[theorem]{Proposition}
\newtheorem{cor}[theorem]{Corollary}

\author{David Callan}
\title{Permutations avoiding 4321 and 3241 have an algebraic generating function}
\affiliation{Department of Statistics,  University of Wisconsin-Madison, USA}
\received{2019-03-15}
\revised{2022-07-25}
\accepted{2022-10-27}

\begin{document}
\publicationdetails{22}{2022}{2}{12}{5286} 
\maketitle
\begin{abstract}
~

We show that permutations avoiding both of the (classical) patterns 4321 and 3241 have the algebraic \gf conjectured by Vladimir Kruchinin.  
\end{abstract}

\section{Introduction} 

This paper is a companion to \cite{kotcallan}, which established the algebraic generating function for  $\{1243,\,2134\}$-avoiding permutations conjectured by Vaclav Kotesovec \cite[\htmladdnormallink{A164651}{http://oeis.org/A164651}]{oeis}. 
In similar vein, Vladimir Kruchinin  \cite[\htmladdnormallink{A165543}{http://oeis.org/A165543}]{oeis} conjectured the \gf
\[
\frac{1}{1 - x\,C\big(x C(x)\big)}
\]
for $\{4321,\,3241\}$-avoiding permutations, where $C(x) := \frac{1-\sqrt{1-4x}}{2x}$ denotes the \gf for the Catalan numbers. 

We will show that \avs do indeed have this \gf. 
First, we use the combinatorial interpretation of the Invert transform to reduce the problem to counting 
\emph{indecomposable}, also known as \emph{sum-indecomposable}, \avs. Then we exhibit a bijective mapping from the set of indecomposable \avs of length $n$ 
to the union of Cartesian products $\bigcup_{k=0}^{n-2}\i_{n-k}(321) \times \c_{k}^{(n-2-k)}$, where $\i_{r}(321)$ is the set of
indecomposable 321-avoiding permutations of length $r$  and $\c_{k}^{(r)}$  is the set of integer sequences 
$(a_{1},a_{2},\dots,a_{k})$ satisfying $1 \le a_{1}\le r+1$ and $1 \le a_{i}\le a_{i-1}+1$ for $i \ge 2$. 
The counting sequences for the sets $\i_{r}$ and $\c_{k}^{(r)}$ are known, and the result follows readily.

Section 2 recalls the notion of indecomposability and the application of the Invert transform to indecomposable
permutations. Section 3 reviews  nonnegative lattice paths and integer sequences whose successive entries 
increase by at most 1, that is, elements of $\c_{k}^{(r)}$. Section 4 defines some notions relevant for our bijection. 
Section 5 presents the main 
bijection and Section 6 explains why it works. Section 7 ties everything together.

\section{Indecomposability and the Invert transform} \label{indecInvert}

A \emph{standard} permutation is one on an initial segment of the positive integers and to \emph{standardize} a
permutation on a set of positive integers means to replace its smallest entry by 1, next smallest by 2 and so on, 
thereby obtaining a standard permutation.  In the context of pattern avoidance, we consider standard 
permutations written in one-line form (that is, as lists). When a standard permutation is written in 
two-line form, it may be possible to insert some vertical bars to obtain subpermutations, not necessarily 
standard, as in 
$\left(\begin{smallmatrix} 1&2&3&|&4&5&6&7 \\
3&1&2&|&6&5&4&7
\end{smallmatrix}\right)$.
After inserting the largest possible number (0 or more) of such bars (to produce nonempty permutations), as in 
$\left(\begin{smallmatrix} 1&2&3&|&4&5&6&|&7 \\
3&1&2&|&6&5&4&|&7
\end{smallmatrix}\right)$,
we obtain the \emph{components} of the permutation, here 3\,1\,2, 6\,5\,4, 7.
A permutation is \emph{indecomposable}, also called \emph{sum-indecomposable} \cite{bevan}, if it has exactly one component. (Thus the permutation 1 is indecomposable but the empty permutation is not.)

Let $F(x) = 1+x+2x^{2}+6x^{3}+\cdots$ denote the \gf for \avs and $G(x)=x+x^{2}+3x^{3}+\cdots$ the \gf for indecomposable \avs.
Clearly, a permutation avoids $\{4321,\,3241\}$ if and only if each of its components does so. 
Hence, the combinatorial interpretation of the Invert transform (see \cite{kotcallan} or \cite{Analytic_Combinatorics}) implies that 
\[
F(x)=\frac{1}{1-G(x)}\, ,
\]
and our problem is reduced to showing that $G(x)=x\,C\big(x C(x)\big)$.

\section{Nonnegative lattice paths} 

It is well known that the ``ballot number'' $\cc{n}{m} := $  
$ \frac{m+1}{2n+m+1}\binom{2n+m+1}{n}$
(with $\cc{n}{-1}:=1 $ if $n=0$ and $:=0$ if $n\ge 1$)  \cite[\htmladdnormallink{A009766}{http://oeis.org/A009766}]{oeis} counts nonnegative paths of
$n+m$ upsteps $U=(1,1)$ and $n$ downsteps $D=(1,-1)$, where nonnegative means the path 
never dips below \gl, the horizontal line through its initial vertex (see, e.g., \cite{cattransform}).
A nonnegative path of $n$ upsteps and $n$ downsteps is a \emph{Dyck} path and its \emph{size} is $n$. A nonempty 
Dyck path is \emph{indecomposable} if its only return to \gl is at the end. The returns to \gl split a Dyck path into its 
(indecomposable) components. The number of indecomposable Dyck paths of size $n$ is $C_{n-1}$ (delete the first and last steps to obtain a one-size-smaller Dyck path).

\begin{prop}$\cite[2.126]{cat2015}$ \label{count1}
The number of indecomposable $321$-avoiding permutations on $[n]$ is $C_{n-1}$.
\end{prop}

\begin{proof} 
One method is to observe that Krattenthaler's bijection \cite{kratt2001} from 321-avoiding 
permutations on $[n]$ to Dyck paths of size $n$ preserves components in the obvious sense and so sends indecomposable permutations
to indecomposable paths. 
\end{proof}

Given a nonnegative path, successively delete the first peak ($UD$) recording its height above \gl 
until no peaks remain, as in Figure~\ref{fig:1} (heights prepended to existing list).

\begin{figure}[htbp]
  \begin{center}
\Einheit=0.4cm
\[
\Pfad(-12,0),3343334344\endPfad
\Pfad(1,0),33334344\endPfad
\SPfad(-12,0),1111111111\endSPfad
\SPfad(1,0),11111111\endSPfad
\Label\o{\rightarrow}(-0.5,2)
\Label\o{\rightarrow}(10.5,2)
\red{
\Pfad(-10,0),22\endPfad
\Pfad(5,0),2222\endPfad}
\DuennPunkt(-12,0)
\DuennPunkt(-11,1)
\DuennPunkt(-10,2)
\DuennPunkt(-9,1)
\DuennPunkt(-8,2)
\DuennPunkt(-7,3)
\DuennPunkt(-6,4)
\DuennPunkt(-5,3)
\DuennPunkt(-4,4)
\DuennPunkt(-3,3)
\DuennPunkt(-2,2)
\DuennPunkt(1,0)
\DuennPunkt(2,1)
\DuennPunkt(3,2)
\DuennPunkt(4,3)
\DuennPunkt(5,4)
\DuennPunkt(6,3)
\DuennPunkt(7,4)
\DuennPunkt(8,3)
\DuennPunkt(9,2)
\Label\o{ \textrm{2}}(-7,-1.6)
\Label\o{ \textrm{$4\:2$}}(5,-1.6)
\]
\\
\[
\Pfad(-9,0),333344\endPfad
\Pfad(0,0),3334\endPfad
\Pfad(7,0),33\endPfad
\SPfad(-9,0),111111\endSPfad
\SPfad(0,0),1111\endSPfad
\SPfad(7,0),11\endSPfad
\Label\o{\rightarrow}(-1.5,1.5)
\Label\o{\rightarrow}(5.5,1)
\red{
\Pfad(-5,0),2222\endPfad
\Pfad(3,0),222\endPfad}
\DuennPunkt(-9,0)
\DuennPunkt(-8,1)
\DuennPunkt(-7,2)
\DuennPunkt(-6,3)
\DuennPunkt(-5,4)
\DuennPunkt(-4,3)
\DuennPunkt(-3,2)
\DuennPunkt(0,0)
\DuennPunkt(1,1)
\DuennPunkt(2,2)
\DuennPunkt(3,3)
\DuennPunkt(4,2)
\DuennPunkt(7,0)
\DuennPunkt(8,1)
\DuennPunkt(9,2)
\Label\o{ \textrm{$4\:4\:2$}}(-6,-1.6)
\Label\o{ \textrm{\framebox{$3\:4\:4\:2$}}}(2,-1.6)
\]
    \caption{Deleting peaks from a Dyck path.}
    \label{fig:1}
  \end{center}
\end{figure}

The path on the left produces the list of heights 3442 and this is a map from nonnegative paths of $n+m$ $U$s and $n$ $D$s to $\c_{n}^{(m)}$.

To reverse the map, suppose given $(a_{1},a_{2}, \dots, a_{n}) \in \c_{n}^{(m)}$. Start with a path of $m\ U$s. Then 
successively insert a peak at height $a_{i}$ into the initial ascent of the current path so that its top vertex 
is at height $a_{i},\ 1\le i \le n$. This produces a path whose first peak is at height $a_{i}$ and 
the growth condition $a_{i+1}\le a_{i}+1$ is just what is needed to enable the next step. So the map is a bijection and we have 
\begin{prop}$\cite{cat2015}$ \label{count2}
 $\mid\! \c_{n}^{(m)} \!\mid\, = \cc{n}{m}$\,.
\end{prop}

The preceding construction will be mirrored in Section \ref{whyworks} below when we insert a ``peak'' entry into a permutation so that 
if it has ``height'' $h$, then there are $h+1$ possibilities for the next insertion.

\section{Some preliminary definitions} 

By (slight) abuse of language, to delete an entry $y \in [n]$ from a permutation $p$ on $[n]$ 
means to erase $y$ and then 
subtract 1 from each entry $>y$; $p\,\textrm{\footnotesize{$\backslash$}}\,\{y\}$ denotes the resulting 
permutation. The non-$y$ entries of $p$ correspond in an obvious way to the entries of 
$p\,\textrm{\footnotesize{$\backslash$}}\,\{y\}$.
Conversely, to insert $y$ in position $i$ means to increment by 1 each entry $\ge y$ 
and then place $y$ in position $i$; we use $p\oplus_{i}y$ to denote the result. Thus, for $p=3142,\,i=4,\,y=2,\ 
p\oplus_{i}\!y=41523$. Again, 
the entries of $p$ correspond naturally to the non-$y$ entries of $p\oplus_{i}\!y$. 
The adjective/noun \lr is short for left-to-right maximum in a permutation.

\textbf{\emph{Henceforth, for brevity, we use the unadorned term ``avoider'' to mean an 
indecomposable $\{\mbf{3241,}$ $\mbf{4321}\}$-avoider.}}

A \emph{key-2} entry in an avoider is an entry that serves as the ``2'' in either a 321 pattern or a 4312 pattern.
For example, the key-2 entries in 6174235 are 4 and 3. Clearly, an avoider with no key-2 entries is a 321-avoider.
The term key-2 is mnemonic but somewhat ungainly and to add a little color, \textbf{we will refer to a key-2 entry as a \emph{blue} entry}.

\begin{figure}[htb]
\vspace{-.1in}
\begin{center}
\includegraphics[angle=0, scale = .8]{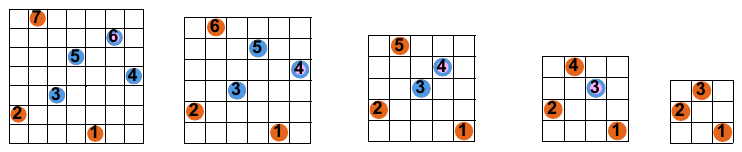}
\end{center}
\vspace{-.2in}
\caption{Avoiders with colored entries and highlight.}
\label{fig:2}
\end{figure}
The \emph{peak} \k entry in an avoider that contains (one or more) 321 patterns is the larger of the last ``1'' of a 321 and its 
(immediate) predecessor. The terminology is justified  because the peak \k is indeed a \k: suppose $a$ 
is the last ``1'' of a 321 in an avoider, 
say the last ``1'' of $cba$, and $y$ is the predecessor of $a$. If $y<a$, then $a$ is the ``2'' of the 4312 pattern  $cbya$. If $y>a$, 
then $y=b$ or $y\ne b$. In the former case, $y$ is the ``2'' of the 321 pattern $cba$; in the latter case, 
$y<c$ (else $cbya$ is a proscribed 3241) and so $y$ is the ``2'' of the 321 pattern $cya$. See Figure~\ref{fig:2} for some examples 
of avoiders with blue entries so colored and peak blue entry highlighted.

We let $\i_{n}(3241,4321)$ denote the set of avoiders (indecomposable $\{3241,4321\}$-avoiding permutations)  of length $n$, and similarly 
$\i_{n}(321)$ is the set of indecomposable $321$-avoiding permutations. 
Set $\i_{n,k}(3241,4321)=\{p \in \i_{n}(3241,4321): p \textrm{ has $k$ \k entries}\}$. 

\section{The bijection} 

\begin{theorem}
For $0\le k \le n-2$,  there is a bijection
\[
\i_{n,k}(3241,4321) \longrightarrow \i_{n-k}(321) \times \c_{k}^{(n-2-k)}\, .
\]
\end{theorem}

\noindent  Here is its description. Suppose given $p \in \i_{n,k}(3241,4321)$. If $k=0$, then $p$ is 
already 321-avoiding and $(p,\eps)$ is the image pair, where $\eps$ denotes the empty list. If $k\ge 1$, 
the idea is to successively delete the (current) peak \k entry recording, in the same right to left fashion as 
in Figure~\ref{fig:1}, its ``height'', appropriately defined, until a 321-avoider $q$ is obtained. Then the 
image pair is $(q,L)$, where $L$ is the list of heights. The trick is to find the correct definition 
of height, and it's a doozy.

To this end, associate to each 321-containing permutation $p$                                                                                                                                                                                                                                                                                                                                                                                                                                                                                                                                                                                                                                                                                                                                                                                                                        
a triple $a<b<c$, all integers except that 
$c$ may be infinite: $a$ is the last ``1'' of a 321 in $p$, $b$ is the 
rightmost entry to the left of $a$ 
that exceeds $a$, and $c$ is the first non-LRMax entry after $a$ (with $c:=\infty$ 
if there is no such non-LRMax). Thus, for $p=321$, we have $(a,b,c)=(1,2,\infty)$ 
and for $p=4631275$, we have $(a,b,c)=(2,3,5)$.

\begin{prop}\label{w}
If $p$ is a $321$-containing avoider with associated triple $(a,b,c)$, 
then there is an entry $w$ in $p$ such that $wba$ is a $321$ pattern in $p$.
\end{prop}

\begin{proof}
Since $a$ is the ``1'' of a 321, there is a sublist $vua$ in $p$ with $v>u>a$. 
By definition of  $b$, $u$ must lie weakly to the left of $b$ and so $v\ne b$. If $v>b$, take $w=v$. 
Otherwise, $v$ and $u$ must both be $<b$ and $vuba$ is a forbidden 3241 pattern. 
\end{proof}

\begin{cor}
If $p$ is a $321$-containing avoider with associated triple $(a,b,c)$, 
then $c>b$.
\end{cor}

\begin{proof}
If not, $wbc$ would be a 321, violating the definition of $a$ as the last ``1'' of a 321. 
\end{proof}

We now define what we call,
for reasons to become clear, the \emph{peak-insertion set} of an avoider.
For a $321$-containing avoider on $[n]$ with associated triple $(a,b,c)$, the peak-insertion set is the disjoint union $[a+1,b+1] \cup [c+1,n]$  
where $[c+1,n]=\emptyset$ if $c=\infty$. For a 321-avoiding avoider on $[n]$,
the $a$ and $b$ evaporate and the peak-insertion set is  
$[c+1,n]$ with $c:=1$, that is, $[2,n]$. Thus, for $p=4631275$ with $(a,b,c)=(2,3,5)$,
the peak-insertion set is $\{3,4,6,7\}$.

Next, we arrange the peak-insertion set of an avoider $p$ into a suitably ordered list, called the \emph{peak-insertion list} of $p$. Taken left to right, the \lr 
entries $> c$ of $p$ form a list $A$ and the non-\lr entries $\ge c$ form a list $B$. Thus $A \cup B = [c,n]$.
Obviously, $A$ is an increasing list, and so is $B$ for otherwise, in the 321-containing case, $a$ would not 
be the last ``1'' of a 321, and in the 321-avoiding case, a 321 would actually
be present. Split $A$ into maximal
runs of consecutive integers $A_{1},A_{2}, \dots, A_{t}$.
Likewise, split $B$ into maximal runs of consecutive
integers but this time written as $b_{1}B_{1},b_{2}B_{2}, \dots , b_{t}B_{t}$, where $b_{i}$ is the first 
entry of the $i$-th run and $B_{i}$ may be empty. There is the same number of 
runs in $A$ as in $B$ because (i) the smallest run contains $b_{1}=c$ and comes from $B$ since $c$ is not a \lr, 
(ii) thereafter the runs alternate between $A$ and $B$, and (iii) the largest run contains $n$, a \lr, and so comes  from $A$.

\begin{figure}
\vspace{-.6in}
\begin{center}
\includegraphics[angle=0, scale = .8]{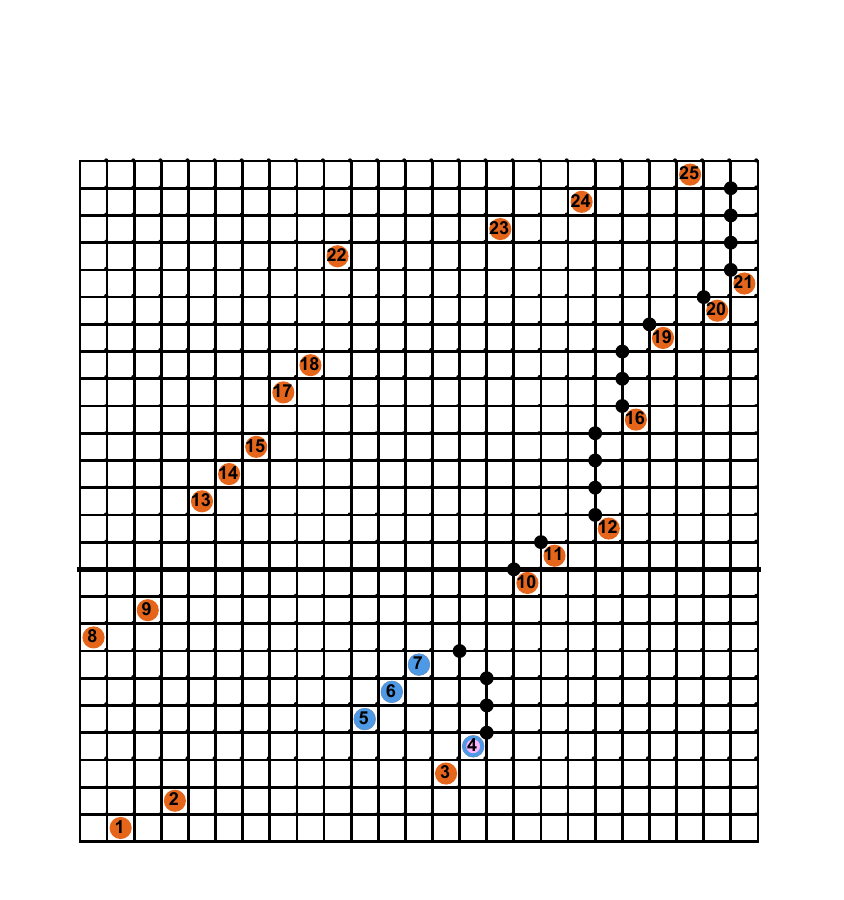}
\end{center}
\vspace{-.4in}
\caption{Avoider with $(a,b,c)=(4,7,10)$ and bullets to illustrate the peak-insertion list.}
\label{fig:3}
\end{figure}

The \emph{peak-insertion list} of $p$ is now defined to be the peak-insertion set of $p$ 
listed in the following order ($L^{r}$ denotes the reversal of the list $L$):
\[
A_{t}\ B_{t}^{r}\ A_{t-1}\ b_{t}\ B_{t-1}^{r}\ A_{t-2}\ b_{t-1}\ B_{t-2}^{r}\  \cdots \  \ A_{1}\ b_{2}\ B_{1}^{r}\ 
b\ \overline{b-1}\ \overline{b-2}\ \cdots\ \overline{a+2}\ \overline{a+1} \ \overline{b+1}\ ,
\]
where the terminal segment starting at $b$ is omitted if $p$ is 321-avoiding.
Note that $b_{1}=c$ is missing and the list consists of $[a+1,b+1] \cup [c+1,n]$\ \ (or $[c+1,n]$ in the 321-avoiding case), as it should.

For example, for the avoider shown in matrix form in Figure~\ref{fig:3}, we have $(a,b,c)=(4,7,10)$ 
and runs in $A$ and $B$ as follows.
\[
\begin{array}{ccccc}
i & = & 1 & 2 & 3 \\[1mm]
A_{i} & = & 13 \ \,14 \ \, 15 & \ \,17\ \, 18 &  \ \,22 \ \,23  \ \, 24  \ \, 25  \\
b_{i}|B_{i} & =  & 10 \ | \ 11 \ \,12  & \ \,16 \ |\ \eps  & \ \,19 \ | \  20 \ \,21
\end{array}
\]
Here, $t=3$ and the ordering in the peak-insertion list is 
\[
\begin{array}{c|c|c|c|c|c|c|c|c|c}
A_{3}& B_3^r & A_2 & b_3 & B_2^r
 & A_1 & b_2 & B_1^r & b \cdots a\!+\!1 & b\!+\!1 \\ \hline
22\ \,23\ \,24\ \,25 & 21 \ \, 20 & 17 \ \,18 &  19 & & 13\ \,14\ \,15 & 16 & 12\ \, 11 & 7\ \,6\ \,5 & 8 
\end{array}\, .
\]

We can now define the \emph{height} of the peak \k entry $y$ in a 321-containing avoider $p$\,: 
it is the position of $y$ in the peak-insertion list of $p\,\textrm{\footnotesize{$\backslash$}}\,\{y\}$. 
(We will see later that $y$ must be in the peak-insertion set of $p\,\textrm{\footnotesize{$\backslash$}}\,\{y\}$.)

\begin{figure}
\vspace{-.2in}
\begin{center}
\includegraphics[angle=0, scale = .8]{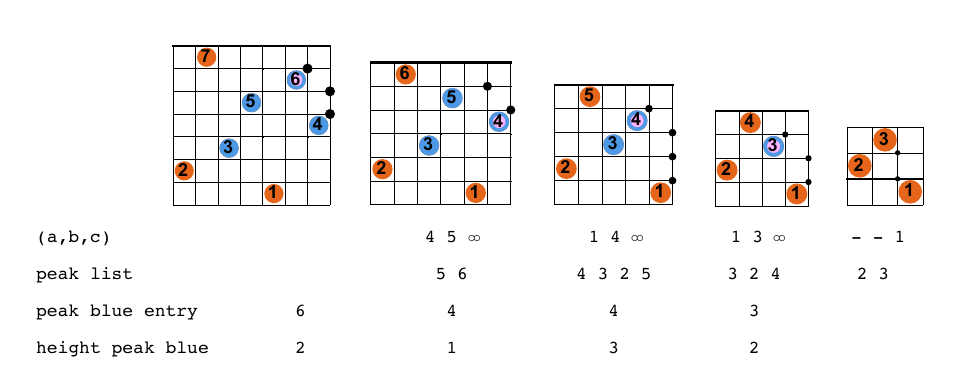}
\end{center}
\vspace{-.2in}
\caption{An example of the bijection.}
\label{fig:4}
\end{figure}
There is a graphical way to visualize the ordering in the peak-insertion list. As illustrated in Figure~\ref{fig:3},
for each $y$ in the peak-insertion set $[a+1,b+1] \cup [c+1,n]$ of $p$, insert a bullet at vertex $(x-1,y-1)$ 
in the matrix diagram of $p$ where the abscissa $x$ is determined as follows. 
For $y>c$, insert the bullet as far right as possible so that the region below and to the right of the bullet is nonempty.
For $y \in [a+1,b+1]$, let $i$ denote the position of $a$ in $p$. Then, 
for $y \in [a+1,b],\ x=i+1$, and for $y=b+1,\ x=i$.
If the bullets are arranged 
in order of distance from the vertical line $x=n$ and, for bullets at the same distance from $x=n$, in order of 
distance from the horizontal line $y=c$ or $y=n$ if $c=\infty$ (heavy line in Figure~\ref{fig:3}), then their $y$'s form the peak-insertion list.

The full mapping is illustrated in Figure~\ref{fig:4}, which shows that $ 2735164 \rightarrow (231,\ 2312)$.

\section{Why it works} \label{whyworks}

We need to establish several facts to show the map does all it claims to and is invertible.
\begin{prop}
Suppose $y$ is the peak \k entry of a $321$-containing avoider $p$.
Then $($i\,$)$ $p\,\textrm{\footnotesize{$\backslash$}}\,\{y\}$ is again an avoider, $($ii\,$)$ the blue entries of 
$p$ other than $y$ become the \k entries of $p\,\textrm{\footnotesize{$\backslash$}}\,\{y\}$, and $($iii\,$)$ $y$ is in the 
peak-insertion set of $p\,\textrm{\footnotesize{$\backslash$}}\,\{y\}$.
\end{prop}

\begin{proof}
(i) $p\,\textrm{\footnotesize{$\backslash$}}\,\{y\}$ inherits the $\{4321,\,3241\}$-avoiding property 
from $p$. If $p\,\textrm{\footnotesize{$\backslash$}}\,\{y\}$ was decomposable then the entries other than $y$ 
in the 321 or 4312 pattern containing $y$ in $p$ would correspond to entries in the same component of 
$p\,\textrm{\footnotesize{$\backslash$}}\,\{y\}$. But then $p$ would also be decomposable, obviously in the 321 case, 
and because the ``1'' and ``2'' can be chosen adjacent in the 4312 case. 
(ii) No new blue entry can be introduced and no non-peak 
blue entry is lost because if the deleted entry $y$ is the ``1'' of a 321, then the ``2'', a blue 
entry in $p$, is still the ``2'' of a 321 in $p\,\textrm{\footnotesize{$\backslash$}}\,\{y\}$ since the predecessor of 
$y$ in $p$  is $< y$ and so serves as a ``1'' in place of $y$. Also, the peak blue entry cannot 
possibly be the ``4'', ``3'', or ``1'' of a 4312, so no blue entry in $p$ that is the ``2'' of a 
4312 loses its blue status in $p\,\textrm{\footnotesize{$\backslash$}}\,\{y\}$. 
(iii) This will be proved in contrapositive form in Proposition \ref{bigone} below. 
\end{proof}

\begin{lemma}\label{misc}
Suppose $p$ is a $321$-containing avoider with associated triple $(a,b,c)$. 

$($i$\,)$ If $c$ is finite, then all entries after $a$ in $p$ are $\ge c$. 

$($ii$\,)$ If $c=\infty$ or $c$ is finite and $c>b+1$, then  $b+1$ lies to the left of $b$ in $p$. 

$($iii$\,)$ If $c=\infty$, then $a$ is the last entry of $p$.

$($iv$\,)$ Suppose $z>b$ is an entry of $p$. Then $z$ is a \lr in $p$ provided $z$ lies to the left of $c$ in $p$ in case $c$ is finite.

\end{lemma}

\begin{proof}
(i) All entries after $c$ are $>c$ else $a$ would not be the last ``1'' of a 321.
If the assertion fails, take $y$ to be the rightmost offending entry in $p$. Clearly, $y$ lies between $a$ and $c$ in $p$ 
and $a<y<c$ and region $Q$ in the schematic of Figure~\ref{fig:5} is empty because $y$ is the rightmost offender.

\begin{figure}[htb]
\vspace{-.1in}
\begin{center}
\includegraphics[angle=0, scale = .95]{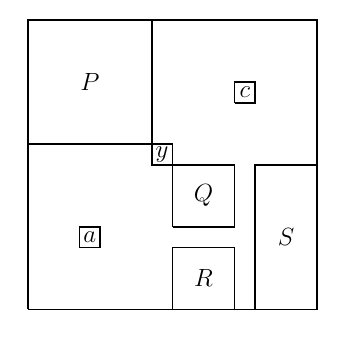}
\end{center}
\caption{Schematic representation of avoider.}
\label{fig:5}
\end{figure}

Also, $P$ is empty since $y$ is a \lr (by definition of $c$), $R$ is empty 
else $a$ is the ``2'' of a 321, and $S$ is empty 
else $c$ is the ``2'' of a 321. These empty regions force $y$ to be a fixed point and $p$ to be decomposable. 

(ii) First, $b+1$ cannot lie between $b$ and $a$ in $p$ by definition of $b$. If $c=\infty$, we are done by part (i).
So suppose $c$ is finite.
The entry $w>b$ whose existence is guaranteed by  Prop.\,\ref{w} implies that $b+1$ is not a \lr, 
violating the definition of $c$ if $b+1$ lies between $a$ and $c$. If $b+1$ lies to the right of $c$, then $c\:
\overline{b+1}$ is the ``21'' of a 321 (since $c$ is not a \lr), 
contradicting the assumption that $a$ is the last ``1'' of a 321.

(iii) If  not, then all entries after $a$ would be \lr entries, and the last entry would be $n$,
violating indecomposability.

(iv) Suppose $z>b$ is an offender. If $c=\infty$, $z$ lies to the left of $b$ by part (iii) and the definition of $b$. 
If $c$ is finite, $z$ lies to the left of $a$ by definition of $c$, and so lies to the left of $b$ by definition of $b$.
In either case, $z$ is a non-\lr lying to the left of $b$. Then $zba$ is the ``321'' of a forbidden 4321. 
\end{proof}

\begin{prop}
For each $y$ in the peak-insertion set of an avoider $p$ on $[n]$, 
there is exactly one position $i$  such that 
$q:=p\oplus_{i}y$ $($insertion of $y$ at position $i$\,$)$ satisfies $($i\,$)$ $q$ is an avoider, 
$($ii\,$)$ the peak \k entry of $q$ is $y$, and 
$($iii\,$)$ $q$ has just one more blue entry than $p$. 
Also, for $y$ not in the peak-insertion set of $p$, there is no such $i$.  
\end{prop}

\begin{proof}
First, suppose $y \in [c+1,n]$. 
Let $z$ be the rightmost entry of $p$ that is $<y$. 
Insert $y$ immediately to the left of $z$. Suppose $p$ has the matrix form depicted schematically 
in Figure~\ref{fig:6} where the bullet represents the inserted entry and $z$ its successor.

\begin{figure}[htb]
\vspace{-.0in}
\begin{center}
\includegraphics[angle=0, scale = .95]{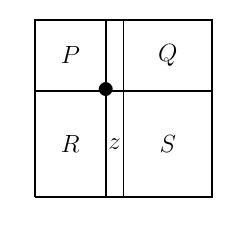}
\end{center}
\vspace{-.2in}
\caption{Schematic for peak-insertion set.}
\label{fig:6}
\end{figure}  
Then $S$ is empty (contains no entries)  by definition of $y$. If $P$ were also empty, $p$ would be decomposable. 
Thus $y$ is the ``2'' of a 321, making $y$ blue in $p\oplus y$ and, clearly, it is the peak blue entry. 
On the other hand, if $y$ is inserted to the right of $z$ it will not be blue, and if inserted to the left of $z$ but not adjacent to $c$, 
it may be blue but will not be the peak blue.

Now suppose $y \in [a+1,b+1]$. If $y=b+1$, insert $y$ just before $a$ ($y$ will be the ``2'' of a 321), 
and if $y \in [a+1,b]$, insert $y$ just after $a$ ($y$ will be the ``2'' of a 4312). Similar considerations show that, 
for this insertion point, $y$ will be the peak blue entry  in $p\oplus y$ and the only new blue entry. 
Also no other insertion point will do.

As for the last assertion, if $p$ is 321-avoiding, the peak-insertion set is $[2,n]$ and 1 cannot be a blue entry in $p\oplus 1$ 
because, by definition, blue entries exceed 1. Now suppose $p$ is 321-containing with associated $a,b,c$. 
If $y\le a$ is inserted to the left of $a$, it cannot be the larger of the last ``1'' of a 321 and its predecessor 
in $p\oplus y$; if inserted to the right of $a$, a descending quadruple is present in $p\oplus y$.
Next, suppose $y\in [b+2,c]$ ($y\ge b+2$ in case $c=\infty$). If $y$ is inserted to the left of $a$, 
then it is not the larger of $a$ and the predecessor of $a$ 
unless it actually is the predecessor of $a$, but in that case $\overline{b+1}\,b\,y\,a$ 
is a forbidden 3241 by Lemma \ref{misc}\,(ii); if $y$ is inserted between $a$ and $c$ in case $c$ is 
finite or after $a$ in case $c=\infty$, 
$y$ cannot be the ``1'' of a 321 in  $p\oplus y$ 
by Lemma \ref{misc}\,(iv) and so $y$ is certainly not the larger of the last ``1'' of a 321 and its 
predecessor in $p\oplus y$; if $y$ is inserted after $c$, then it is the last ``1'' of a 321 in $p\oplus y$, 
but is not larger than its predecessor and so is not peak.  
\end{proof}

\begin{prop}\label{bigone}
For an avoider $p$, as $y$ ranges from left to right over the peak-insertion list of $p$, 
the length of the 
peak-insertion list of $p\oplus_{i_{y}}\!y$\, ranges from left to right over the interval $2,3, \dots, r+1$, 
where $r$ denotes the length of the peak-insertion list of $p$ and $i_y$ is the $i$ of the preceding Proposition.
\end{prop}

\begin{proof}
Recall that every avoider $p$ is associated with an $(a,b,c)$ triple if it is 321-containing and 
with a singleton $c=1$ otherwise, and the   
peak-insertion set for $p$ is $[a+1,b+1] \cup [c+1,n]$
with the first interval absent if $p$ is 321-avoiding and the second interval absent if $c=\infty$. 
We need to determine the triples, denoted $(a_y,b_y,c_y)$, for each $p\oplus_{i_{y}}\!y$ with $y$ in the peak-insertion list of $p$. 
In the peak-insertion list of $p$, the entries $>c$ all occur before the entries $<c$.
\begin{figure}[htb]
\vspace{-.1in}
\begin{center}
\includegraphics[angle=0, scale = .7]{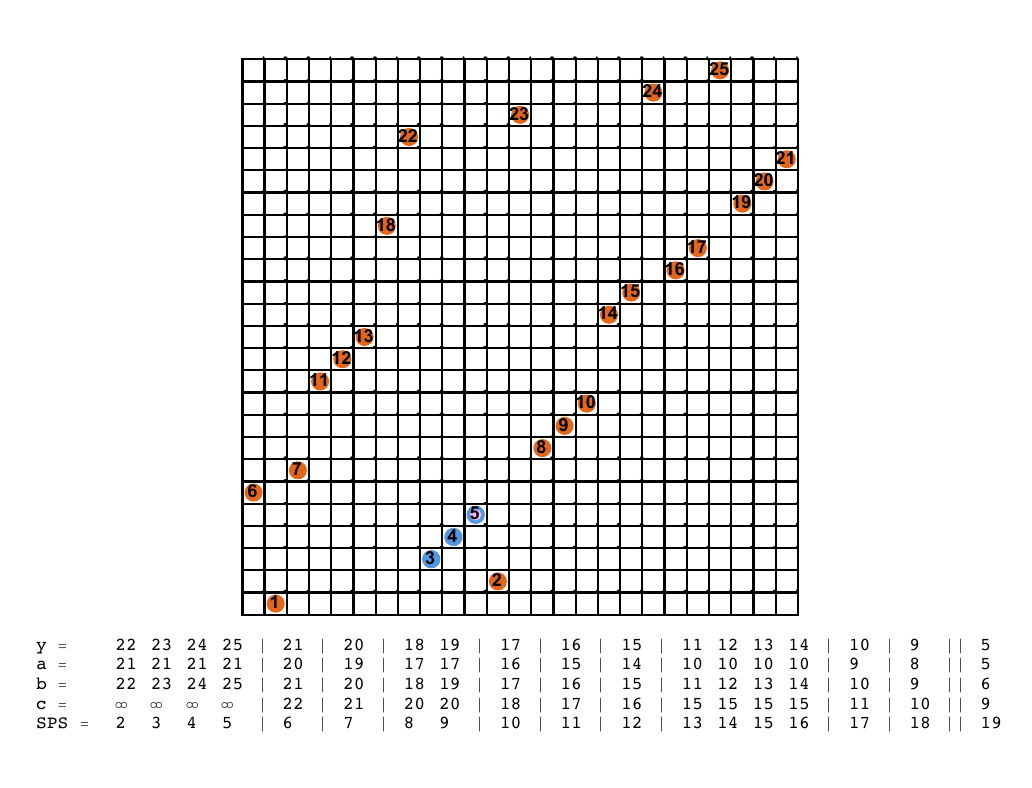}
\end{center}
\vspace{-.3in}
\caption{Example showing $a,b,c$ and size of peak-insertion set for each $y$.}
\label{fig:7}
\end{figure}
Split the entries $>c$ into segments consisting of increasing runs as illustrated in Figure~\ref{fig:7}  
for an avoider on $[n]=[25]$ with $(a,b,c)=(4,7,10)$ and peak-insertion set $[a+1,b+1] \cup [c+1,n] = [5,8] \cup [11,25]$. 
Then for $y>c$ in the peak-insertion list of $p$,
the $a_y,b_y$ and $c_y$ of $p\oplus_{i_{y}}\!y$ are as follows: 
(i) $a_y$ is one less than the smallest entry in the segment containing $y$, 
(ii) $b_y=y$, and 
(iii) $c_y$ is one more than the largest entry in the segment containing $y$ (or $\infty$ if this largest entry is $n$).
For $y<c$ in the peak-insertion list of $p$,  
the $b_y$ is $b+1$, the $c_y$ is $c+1$, while $a_y = y$ for $y\in [a+1,b]$ and $a_y=a$ for $y=b$. 
We leave the reader to verify the truth of these assertions with the visual aid that
each pair $(i_y,y)$ is shown as a bullet at vertex $(i_{y}-1,y-1)$ in Figure 7. So, ``expanding'' the $(i_y,y)$ bullet 
into a cell containing the entry $y$ gives the matrix diagram of $p \oplus_{i_{y}} y$.

It is now clear that the size of the peak-insertion set (SPS in Figure 7) starts at 2 
and increases by 1 thereafter as $y$ ranges across the peak-insertion list of $p$. 
\end{proof}

\section{Putting it all together} 

From Propositions \ref{count1} and \ref{count2} and the preceding bijection, we find that the number $u_n$ of indecomposable 
$\{4321,\,3241\}$-avoiding permutations of length $n$ is given by $u_0=0,\,u_1=1$ and, for $n\ge 2$,
\begin{eqnarray*}
u_n & = & \sum_{k=0}^{n-2}C_{n-1-k}C_{k,n-2-k} \\
 & = & \sum_{k=0}^{n-2}\frac{n - 1 - k}{n - 1 + k}\binom{n - 1 + k}{ k} C_{n-1-k}        \, .
\end{eqnarray*}
This is sequence  
\htmladdnormallink{A127632}{http://oeis.org/A127632} in OEIS \cite{oeis} with \gf
\[
\sum_{n\ge 0}u_n x^n = x\,C\big(x C(x)\big),
\]
where $C(x)$ is the \gf for the Catalan numbers, and by Section \ref{indecInvert}, 
the claimed \gf for $\{4321,\,3241\}$-avoiding permutations follows. The counting sequence for 
$\{4321, 3241\}$-avoiders can be succinctly described 
as the Invert transform of the Catalan transform of the Catalan numbers.

\bigskip

{\large \textbf{Addendum}}\quad 
After this paper was posted to the ArXiv in 2013, a more transparent proof was obtained by Bloom and Vatter \cite{vignettes}.
\vspace*{5mm}

\acknowledgements
I thank an anonymous referee for several helpful remarks.


\begin{thebibliography}{99}

\bibitem{kotcallan} David Callan, The number of \{1243,\:2134\}-avoiding permutations, preprint,
\htmladdnormallink{http://front.math.ucdavis.edu/1303.3857}{http://front.math.ucdavis.edu/1303.3857}, 2013.

\bibitem{oeis} The On-Line Encyclopedia of Integer Sequences, published electronically at 
\htmladdnormallink{http://oeis.org}{http://oeis.org}, 2022.

\bibitem{bevan} David Bevan, Permutation patterns: basic definitions and notation, preprint, \htmladdnormallink{https://arxiv.org/abs/1506.06673}{https://arxiv.org/abs/1506.06673}, 2015.

 
\bibitem{Analytic_Combinatorics} Philippe Flajolet and Robert Sedgewick, 
\htmladdnormallink{\emph{Analytic Combinatorics}}{http://algo.inria.fr/flajolet/Publications/AnaCombi/book.pdf}, Cambridge University Press, 2009.

\bibitem{cattransform} 
David Callan, A combinatorial interpretation of the Catalan transform of the Catalan numbers, preprint,
\htmladdnormallink{http://front.math.ucdavis.edu/1111.0996}{http://front.math.ucdavis.edu/1111.0996}, 2011.

\bibitem{cat2015} Richard Stanley, \emph{Catalan Numbers}, Cambridge University Press, 2015.

\bibitem{kratt2001} Christian  Krattenthaler, Permutations with restricted 
patterns and Dyck paths,
\ \textit{Advances in Applied Mathematics}, \textbf{27}, 2001,  no. 2-3, 510--530. 

\bibitem{vignettes} Jonathan Bloom and Vincent Vatter, Two vignettes on full rook placements, \emph{Australasian Journal of Combinatorics}
\textbf{64} (1) (2016), 77--87.


 \end{thebibliography}
\end{document}